\documentclass[10pt,a4paper,final]{amsart}
\usepackage[T1]{fontenc}
\usepackage[utf8]{inputenc}
\usepackage[english]{babel}
\usepackage{mathrsfs}

\usepackage{graphicx}
\usepackage{import}			    
\usepackage{calc}				
\usepackage{tikz-cd}

\usepackage{todonotes}
\usepackage{comment}
\usepackage{xcolor}			    
\usepackage[notcite,notref]{showkeys}

\usepackage{braket}			    
\usepackage{mathtools}			
	\mathtoolsset{centercolon}	

\usepackage[autostyle,english=british]{csquotes}
\usepackage[style=numeric-comp,bibencoding=utf8]{biblatex}
	\renewbibmacro{in:}{\ifentrytype{article}{}{\printtext{\bibstring{in}\intitlepunct}}}
\usepackage[colorlinks]{hyperref}
	\hypersetup{hidelinks}		
\usepackage{bookmark}
\usepackage[capitalise,noabbrev]{cleveref}

\theoremstyle{plain}
\newtheorem{thm}{Theorem}

\newtheorem{lemma}[thm]{Lemma}

\theoremstyle{definition}
\newtheorem{defi}[thm]{Definition}

\theoremstyle{remark}

\numberwithin{equation}{section}

\newcommand{\numberset}{\mathbb}
\newcommand{\CC}{\numberset{C}}

\newcommand{\NN}{\numberset{N}}
\newcommand{\PP}{\numberset{P}}
\newcommand{\QQ}{\numberset{Q}}
\newcommand{\RR}{\numberset{R}}
\renewcommand{\SS}{\numberset{S}}
\newcommand{\TT}{\numberset{T}}
\newcommand{\ZZ}{\numberset{Z}}

\DeclareMathOperator{\SL}{SL}
\DeclareMathOperator{\GL}{GL}

\DeclareMathOperator{\id}{id}

\DeclareMathOperator{\Realpart}{Re}
\renewcommand{\Re}{\Realpart}

\renewcommand{\epsilon}{\varepsilon}
\renewcommand{\phi}{\varphi}

\newcommand{\cH}{\mathcal{H}}
\newcommand{\cO}{\mathcal{O}}

\begin{document}

\title[Octagonal continued fraction and diagonal changes]{Octagonal continued fraction and \\ diagonal changes}

\author[M. Artigiani]{Mauro Artigiani}
\address{Escuela de Ingeniería, Ciencia y Tecnología\\
 				 Universidad del Rosario\\
				 Bogot\'a\\
				 Colombia}
\email{mauro.artigiani@urosario.edu.co}



\date{\today}

\begin{abstract}
	In this short note we show that the octagon Farey map introduced by Smillie and Ulcigrai in~\cite{SmillieUlcigrai:sturmian,SmillieUlcigrai:geodesic} is an acceleration of the diagonal changes algorithm introduced by Delecroix and Ulcigrai in~\cite{DelecroixUlcigrai:diagonalchanges}.
\end{abstract}


\maketitle

\section{Introduction}
The theory of continued fractions is a beautiful page of mathematics which connects number theory, (hyperbolic) geometry and dynamical systems.
Given a number $\alpha\in\RR$, its continued fraction expansion is an expression of the form
\[
	\alpha=[a_0;a_1,a_2,\dots]=a_0+\cfrac{1}{a_1+\cfrac{1}{a_2+\dots}},
\]
where $a_0\in\ZZ$ and $a_i\in\NN$, for $i\neq 0$.
The rational approximations $p_n/q_n=[a_0;a_1,\dots,a_n]$ obtained by truncating the continued fraction at level $n$ are called \emph{convergents} and are the best approximations to the number $\alpha$, among the ones with denominator bounded by $q_n$.
Subtracting the integer part of $\alpha$, we can assume that $\alpha\in [0,1]$.
The continued fraction of $\alpha$ can then be obtained from the itinerary of the Gauss map $G(x)=\bigl\{\frac{1}{x}\bigr\}$ on $[0,1]$, where $\{\cdot\}$ denotes the fractional part
\[
	a_i=n\iff G^{i-1}(x)\in\left(\frac{1}{n+1},\frac{1}{n}\right].
\]

\begin{figure}[tb]
	\centering
	\def\svgscale{0.6}
	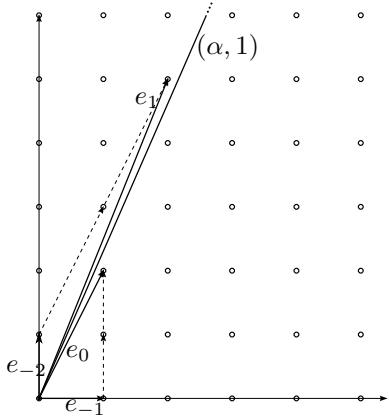
	\caption{The geometric construction of the convergents of $\alpha\in\RR^+$.}
	\label{fig:convergents}
\end{figure}

The continued fraction algorithm can be also realized in a geometric fashion in the following way, see the introduction of~\cite{DelecroixUlcigrai:diagonalchanges} for more details.
Having chosen $\alpha\in\RR^+\setminus\QQ$, we draw the line in direction $(\alpha,1)$.
Then we consider the basis of $\ZZ^2$ given by the vectors $e_{-2}=(0,1)$ and $e_{-1}=(1,0)$.
Note that the line in direction $(\alpha,1)$ is contained in the cone generated by the vectors $e_{-1}$ and $e_{-2}$.
At each step $n\geq 0$, we are going to replace $e_{n-2}$ with a new vector $e_n$ obtained by adding to the vector $e_{n-2}$ the vector $e_{n-1}$ as many times as we can \emph{without} crossing the line in direction $(\alpha,1)$, see \cref{fig:convergents}.
In other words
\[
	e_n=a_n e_{n-1}+e_{n-2}.
\]
This shows that, after the step $n=1$ when we have replaced both our starting vectors, the algorithm is selecting the points in the integer lattice $\ZZ^2$ that are the closest ones to the line $(\alpha,1)$ up to their given height.
Moreover, it follows from the construction that at each step $e_n$ and $e_{n-1}$ form a basis of $\ZZ^2$ and that the line in direction $(\alpha,1)$ is contained in the cone generated by them.

One can show that this procedure produces the continued fraction of $\alpha=[a_0;a_1,\dots]$ and that, if $e_n=(p_n,q_n)$, then $p_n/q_n$ is the $n^{\text{th}}$ convergent to $\alpha$.

It is worth to mention that intermediate vectors of the form $i e_{n-1}+e_{n-2}$, for $i=1,\dots,a_n-1$ are also of interest.
In fact they yield the additive continued fraction convergents, that is the ones produced by the Farey map, whose acceleration gives the Gauss map itself.
These intermediate convergents are called \emph{approximations of the first kind} in the literature, see~\cite{Khinchin:CF}.

It is well-known that the classical continued fraction algorithm acts as a \emph{renormalizing operator} on irrational rotations of the unit interval.
It is easy to see that the induced transformation on a Poincaré section of the geodesic flow in an irrational direction on the flat torus $\TT^2=\RR^2/\ZZ^2$ is an irrational rotation.
Hence, one can use the Gauss map to renormalize the geodesic flow on the flat torus.
From a different point of view, the continued fraction arises from a Poincaré section for the geodesic flow on the moduli space of flat tori, which is (the unit tangle bundle to) a hyperbolic surface, see~\cite{Series:Modular}.

Translation surfaces are higher genus analogues of flat tori, defined by gluing a set of polygons in the plane via translations, see \cref{sec:definitions}.
Translation surfaces carry a Euclidean structure, hence the geodesic flow on any such surface is given, as in the case of the torus, by a straightline flow in a fixed direction.
It is easy to see that the first return map to a transversal for the straightline flow is an interval exchange transformation, which are a generalization of rotations.

It is natural to generalize the theory of continued fractions to translation surfaces.
One way to do this is via Rauzy-Veech induction on interval exchange transformations, see~\cite{Yoccoz:continuedfractions,Zorich:flat}.
Another point of view, which is a direct generalization of the flat geometric point of view on continued fraction described at the beginning of this introduction, has been taken by Delecroix and Ulcigrai in~\cite{DelecroixUlcigrai:diagonalchanges} and their \emph{diagonal changes algorithm} for translation surfaces living in the \emph{hyperelliptic component}.
We will recall the basic definitions of diagonal changes in \cref{sec:diagonalchanges} below.

A particular family of translation surface is the one of \emph{Veech surfaces} (also called lattice surfaces), originally discovered in~\cite{Veech:Eisenstein}.
Examples of Veech surfaces are the surfaces obtained from gluing opposite sides of a regular $2n$-gon in the plane by translation.
By definition, the moduli space of affine deformations of a Veech surface is also (the unit tangle bundle to) a hyperbolic surface.
Hence, one can use methods inspired by hyperbolic geometry, such as the classical ones by Bowen and Series in~\cite{BowenSeries,Series:GeometricalCoding}, to code the geodesic on the moduli space of affine deformations of a Veech surface and deduce a continued fraction algorithm from this construction.

Using this point of view Smillie and Ulcigrai have introduced in~\cite{SmillieUlcigrai:sturmian,SmillieUlcigrai:geodesic} a continued fraction algorithm for the translation surface obtained from the regular octagon (and more generally for all regular $2n$-gons).
Their algorithm can be used to study the straightline flow on the regular octagon from a symbolic point of view, and comes from a particular section of the geodesic flow on the moduli space of affine deformations of the regular octagon.
A nice feature of their algorithm is that, unlike the ones defined by Bowen and Series, behaves as a \emph{full shift} on $7$ symbols, apart from the first move.

On the surface obtained by gluing opposite sides of a regular octagon in the plane by translation, both the diagonal changes algorithm and the Smillie-Ulcigrai algorithm can be used.
Since they are both generalization of the classical continued fraction algorithm on the torus, it is natural to ask whether they are related or not.
The content of this note is to show that indeed they are.

\begin{thm}\label{thm:main}
	The octagon additive continued fraction algorithm defined in~\cite{SmillieUlcigrai:sturmian} is an acceleration of the diagonal changes algorithm for the octagon itself.
\end{thm}

Since, as we remarked above, the continued fraction on the octagon is morally a full-shift, this result allows, to a great extent, to bypass the combinatorial complexity of the diagonal changes algorithm, restricting the analysis to a family of loops in the graph of the induction corresponding to the basic moves of the Smillie-Ulcigrai algorithm.

We remark that is an open question to characterize the behavior of diagonal changes on a Veech surface.

\subsection*{Organization of the paper}
In \cref{sec:definitions} we recall the definitions we need about translation surfaces.
Then we proceed to describe the additive continued fraction algorithm defined in~\cite{SmillieUlcigrai:sturmian,SmillieUlcigrai:geodesic}.
In \cref{sec:diagonalchanges} we recall the definitions for diagonal changes, and we give a different combinatorial description for the octagon, which is more suited for our discussions.
Finally in \cref{sec:thm} we show \cref{thm:main}.
The drawings needed are included in an Appendix at the end of the document.

\section{Definitions}\label{sec:definitions}
We now introduce the basic definitions on translations surfaces which will be needed in the next sections.
General reference on the subject are~\cite{Masur:ergodicsurfaces,Zorich:flat,ForniMatheus:Introduction}.

A compact \emph{translation surfaces} is a finite collection of polygons $\{ P_1, \ldots, P_n\}$ embedded in the plane $\RR^2\cong\CC$ together with side identifications as follows.
Every side $s_i\in P_i$ is identified with a unique side $s_j\in P_j$ such that the sides $s_i$ and $s_j$ are parallel and have the same length.
Moreover, the outward pointing normal vectors with respect to the two sides point in opposite directions.
We then identify the sides $s_i$ and $s_j$ by translations.
We denote by $X$ the surface obtained after performing all the gluing.

We remark that the presentation of a translation surface as a collection of polygons is \emph{not} canonical.
In fact, two collections that differ by \emph{cut and paste} yield the same surface.
More precisely, a ``cut'' operation means cutting some polygon(s) along a straight line connecting two vertices, recording in the new collection of polygons that those sides that have been created are identified in the quotient; a ``paste'' operation corresponds to gluing some polygons along sides that are identified in the quotient.
Two translation surfaces $X=\{P_1,\ldots, P_n\}$ and $X'=\{P'_1,\ldots, P'_m\}$ are \emph{isomorphic} if there exists a (finite) sequence of cut and paste operation that transforms the colletion $\{P_1,\ldots, P_n\}$ into $\{P'_1,\ldots, P'_m\}$, with the appropriate side identifications.
Cut and paste operations are at the heart of the diagonal changes algorithm, which will be described in \cref{sec:diagonalchanges}

The surface $X$ inherits everywhere except in a finite set $\mathscr{S}$, which is contained in the image of the vertices of the polygons, the Euclidean structure from $\RR^2$.
These points are called \emph{conical singularities}.
Around a point $s\in\mathscr{S}$ the total angle is $2\pi (k_s +1)$ for $k_s\in\NN$.
One has the following Gauss-Bonnet formula for the flat metric on the surface:
\[
	2g - 2 =\sum_{s\in\mathscr{S}} k_s,
\]
where $g$ is the genus of the surface $X$.

The collection of translation surfaces with the same topology, that is number of singularities and value of conical angle around each of them (and hence the same genus), is called a \emph{stratum} and is denoted $\cH(k_1, \ldots, k_n)$.
One can show that strata are complex orbifold, not necessarily connected.

Thanks to the Euclidean structure on $X$, for every angle $\theta\in\SS^1$ we have a well-defined concept of linear flow in direction $\theta$, which is given in charts by following lines in direction $\theta$ on $X$.
This corresponds to the geodesic flow on $X$.
A trajectory of the linear flow that connects two (not necessarily distinct) singularities and contains no singularities in its interior is called a \emph{saddle connection}.
To a saddle connection we can associate a \emph{displacement vector} (often called \emph{holonomy vector}), by developing the saddle connection to the plane $\RR^2$ and taking the difference of its endpoints.
In the following, for simplicity, we will often identify saddle connections with their respective displacement vector.
A \emph{separatrix} is a trajectory of the linear flow with only one of its endpoints in a singularity.

There is a natural action by affine diffeomorphisms of $\GL(2,\RR)$, of  on translation surfaces, given by acting on the polygons that constitute the surface by linear transformation.
As the action of $\GL(2, \RR)$ preserves parallelism, this descents to an action on the surface itself.
One can show that the action is continuous on each stratum (with respect to the orbifold topology).
The group of affine diffeomorphisms of a translation surface is called the \emph{Veech group} of $X$.
The Veech group is a discrete subgroup of $\SL_\pm(2,\RR)$, the matrices with determinant equal to $\pm 1$.
A surface is called a \emph{Veech surface} if its Veech group is a lattice inside $\SL_\pm(2,\RR)$.
We remark that we will allow orientation reversing affine diffeomorphisms, as this will allow to use the full dihedral group of the regular octagon in \cref{sec:octagonFarey}.

\section{The octagon Farey map}\label{sec:octagonFarey}
In this section we will recall the definition of the octagon Farey map.
Our presentation will closely follow the one given in~\cite{SmillieUlcigrai:geodesic}.

Let $\cO\subset\CC$ be a regular octagon.
We will use $X=X_\cO$ to denote the translation surface obtained by gluing opposite parallel sides of the octagon.
This surface has genus $2$ and a single conical singularity of order $6\pi$, coming from the image of the vertices of $\cO$, hence it belongs to the stratum $\cH(2)$.

Let $D_8\in\GL(2,\RR)$ the dihedral group of $\cO$, that is the full group of symmetries of the regular octagon.
The octagon Farey map will act as a renormalization operator on $\SS^1$, the space of directions of trajectories.
Since $-\id\in D_8$, we can restrict our analysis to the upper half $\Sigma_+$ of $\SS^1$.
More precisely $\Sigma_+$ is the part corresponding to complex numbers with positive imaginary part.
As we are thinking of $\SS^1$ as the space of directions, it is more convenient to use angle coordinates $\theta\in[0,2\pi)$ to parametrize points $z=e^{i \theta}$.
In other words the angle $\theta$ corresponds to the unit vector $(\cos\theta,\sin\theta)$ in $\RR^2$.
In this coordinates, $\Sigma_+$ corresponds to $\theta\in[0,\pi]$.
Another coordinate we are going to use is the inverse slope coordinate $u$ on $\Sigma_+$ given by $u=\cot(\theta)$.
It is natural in this context to extend $u$ to a map from $\Sigma_+$ to $\RR\PP^1=\RR\cup\{\infty\}$ sending the endpoints of $\Sigma_+$ to the point at infinity.
This coordinate is helpful for us since it allows to conveniently express the action of $\GL(2,\RR)$ on $\SS^1$ simply by M\"obius maps in the $u$ coordinate.

We divide $\Sigma_+$ into $8$ sectors $\overline{\Sigma}_j=\Set{\theta\in\SS : \frac{j\pi}{8}\leq\theta\leq\frac{(j+1)\pi}{8}}$, for $i=0,\dots,7$.
The sector $\overline{\Sigma}_0$ is a fundamental domain for the action of $D_8$ on $\Sigma_+$.
We denote by $\nu_j\in D_8$ the element mapping linearly each sector $\overline{\Sigma}_j$ onto $\overline{\Sigma}_0$.
One can check that these elements are
\[
	\begin{aligned}
		\nu_0&=\begin{pmatrix} 1 & 0 \\ 0 & 1\end{pmatrix}, &
		\nu_1&=\begin{pmatrix} \frac{1}{\sqrt{2}} & \frac{1}{\sqrt{2}} \\ \frac{1}{\sqrt{2}} & -\frac{1}{\sqrt{2}}\end{pmatrix}, &
		\nu_2&=\begin{pmatrix} \frac{1}{\sqrt{2}} & \frac{1}{\sqrt{2}} \\ -\frac{1}{\sqrt{2}} & \frac{1}{\sqrt{2}}\end{pmatrix}, &
		\nu_3&=\begin{pmatrix} 0 & 1 \\ 1 & 0\end{pmatrix}, \\
		\nu_4&=\begin{pmatrix} 0 & 1 \\ -1 & 0\end{pmatrix}, &
		\nu_5&=\begin{pmatrix} -\frac{1}{\sqrt{2}} & \frac{1}{\sqrt{2}} \\ \frac{1}{\sqrt{2}} & \frac{1}{\sqrt{2}}\end{pmatrix}, &
		\nu_6&=\begin{pmatrix} -\frac{1}{\sqrt{2}} & \frac{1}{\sqrt{2}} \\ -\frac{1}{\sqrt{2}} & -\frac{1}{\sqrt{2}}\end{pmatrix}, &
		\nu_7&=\begin{pmatrix} -1 & 0 \\ 0 & 1\end{pmatrix}.
	\end{aligned}
\]
Using these maps, we define a folding map $\operatorname{fold}\colon\Sigma_+\to\overline{\Sigma}_0$ that sends a point $\theta\in\overline{\Sigma}_j$ to the point $\nu_j(\theta)$ with the linear action of $\nu_j$ on the corresponding unit vector $(\cos\theta,\sin\theta)$.
The different branches of $\operatorname{fold}$ agree on the common endpoints and hence we see that $\operatorname{fold}$ is a continuous, piecewise linear, map.

Consider now the element
\[
	\gamma=\begin{pmatrix} -1 & 2(1+\sqrt(2))\\ 0 & 1\end{pmatrix}.
\]
One can show that $\gamma$ and $D_8$ generate the whole Veech group of $X$.
We remark that $\gamma^2=\id$.
If we denote with $\overline{\Sigma}=\overline{\Sigma}_1\cup\dots\cup\overline{\Sigma}_7$, we see that $\gamma$ maps $\overline{\Sigma}_0$ to $\overline{\Sigma}$, and vice versa, reversing the orientation.

\begin{figure}
\centering
\includegraphics[height=0.2\textheight]{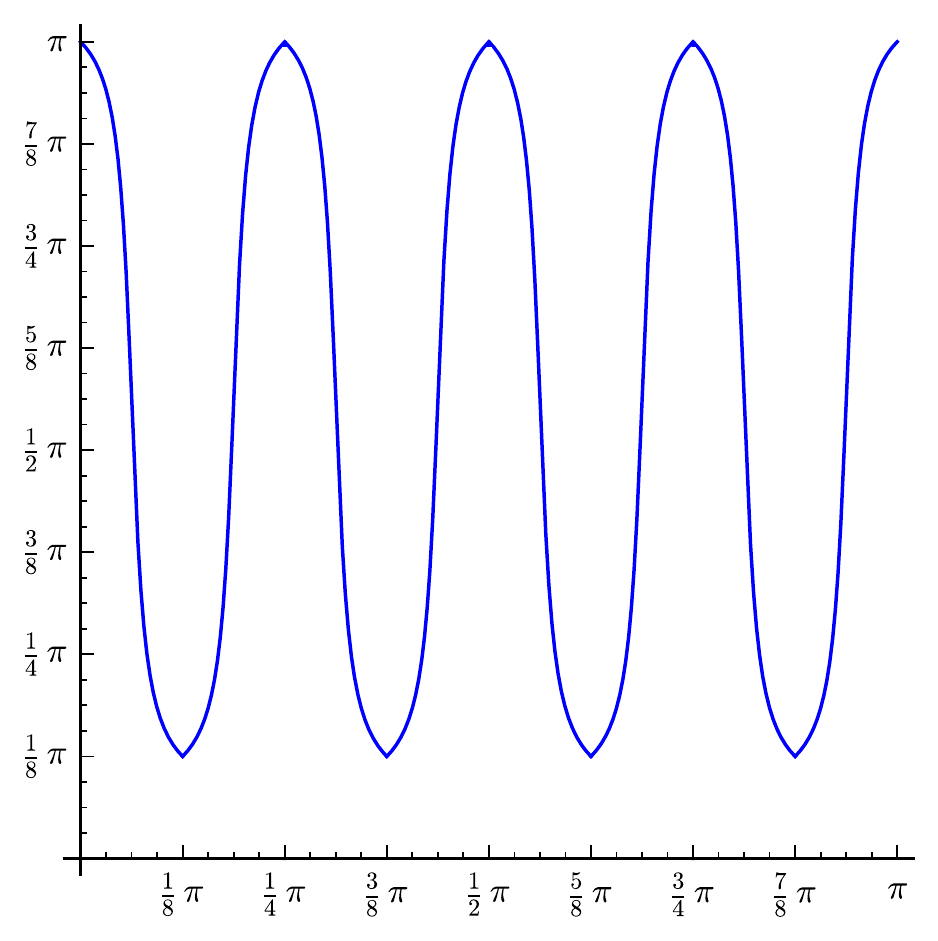}
\caption{The octagon Farey map in angle coordinates.}
\label{fig:Farey}
\end{figure}

Call $F_i\colon\overline{\Sigma}_i\to\RR\PP^1$ the map induced by $\gamma\nu_i$.
We define the \emph{octagon Farey map} $F\colon\RR\PP^1\to\RR\PP^1$ to be the map that acts on directions belonging to the sector $\overline{\Sigma}_i$ as $F_i$.
In other words $F=\gamma\circ\operatorname{fold}$, see \cref{fig:Farey}.
This, in turn, implies that $F$ is a continuous map.
As we said above, the action of $F$ is expressed in the inverse slope coordinate $u$ simply by M\"obius transformation: if $u\in\overline{\Sigma}_i$ we have
\[
	F(u)=\gamma\nu_i*u=\frac{au+b}{cu+d}, \qquad \text{where}\qquad \gamma\nu_i=\begin{pmatrix}
								a & b \\ c & d
							\end{pmatrix}
\]

The action in the angle coordinate is obtained by conjugation with $\cot$.
In the $\theta$ coordinate the map $F$ is expanding at every point, except at the endpoints of each sector, but the amount of expansion is not uniform and tends to one at the endpoints of each sector.
Since all $F_i$ are monotonic, we can define their inverses $F_i^{-1}\colon\overline{\Sigma}\to\overline{\Sigma}_i$, for $i=0,\dots,7$.

We are now ready to recall the definition of an additive continued fraction algorithm, exploiting the map $F$.
Take a direction $\theta\in[0,\pi]$ and record its \emph{itinerary} $\{s_k\}_{k\in\NN}$ under the map $F$.
In other words, we write $s_k=j$ if and only if $F^k(\theta)\in\overline{\Sigma}_j$.
This itinerary is unique if $F^k(\theta)$ never coincides with the endpoint of two sectors.
We remark that, as the image of $F$ is contained in $\overline{\Sigma}$, only $s_0$ can be $0$.
On the other hand, given a sequence $\{s_k\}_{k\in\NN}$ of entries $0,\dots,7$ such that $s_k=0$ implies $k=0$, we consider the intersection $\cap_{k\in\NN} F_{s_0}^{-1}F_{s_1}^{-1}\dots F_{s_k}^{-1}[0,\pi]$.
One can show that the intersection is non empty and consists of only one point.
We hence write
\begin{equation}\label{eq:octagoncf}
	\theta=[s_0;s_1,s_2,\dots]_\cO:=\bigcap_{k\in\NN} F_{s_0}^{-1}F_{s_1}^{-1}\dots F_{s_k}^{-1}[0,\pi],
\end{equation}
for an \emph{octagon Farey expansion of $\theta$}.

One direction $\theta$ can have at most two expansions.
In fact, let us call \emph{terminating} a direction whose continued fraction entries $s_k$ are eventually all $1$ or $7$.
Then, all points that are not endpoints of a sector $\overline{\Sigma}_j$ have a unique expansion.
More precisely, the two sequences $(\ldots,s_k,1,1,1,\ldots)$ and $(\ldots,s_k+1,1,1,1,\ldots)$ correspond to the same direction if $s_k$ is even and  $(\ldots,s_k,7,7,7,\ldots)$ and $(\ldots,s_k+1,7,7,7,\ldots)$ correspond to the same direction if $s_k$ is odd.
Finally $0$ corresponds to $[0;7,7,7,\ldots]_\cO$ and $\pi=[7;7,7,7,\ldots]_\cO$.

Since each $\gamma\nu_i$ maps the corresponding sector $\overline\Sigma_i$ onto $\overline\Sigma$, we want to think of the octagon Farey map $F$ as a \emph{renormalization scheme} acting on directions $\theta\in [0,\pi]$.
To illustrate what we mean by this, let us consider a direction $\theta$ and let us suppose that its first entry in the octagon Farey expansion is not zero.
Then $\theta$ belongs to some $\overline\Sigma_i\subset\overline\Sigma$.
Apply $F$ to $\theta$ hence corresponds to apply the map $F_i=\gamma\nu_i$, which opens up the sector $\overline\Sigma_i$ onto the union of possible sectors $\overline\Sigma$.
By construction, $F(\theta)$ still belongs to $\overline\Sigma$.
Moreover, it is clear from~\eqref{eq:octagoncf} that $F$ acts on the Farey expansion of $\theta$ as a left shift.
In other words, if $\theta=[s_0; s_1, s_2,\ldots]_\cO$ then $F(\theta)=[s_1; s_2, s_3,\ldots]_\cO$.

In the following, given a direction $\theta=[s_0;s_1,s_2,\ldots]_\cO$, we will abuse the notation and continue to call the octagon Farey map the sequence of affine diffeomorphisms given by the octagon continued fraction expansion of $\theta$.

\section{The diagonal changes algorithm}\label{sec:diagonalchanges}
\subsection{Basic definitions}
We are now going to recall the basic definitions of the diagonal changes algorithm, as defined in~\cite{DelecroixUlcigrai:diagonalchanges}.
For more details and for applications of this algorithm we refer the reader to their original paper.

The diagonal changes algorithm produces a sequence of saddle connections which approximate a given direction $\theta\in\SS^1$.
These saddle connections from a \emph{wedge}, in the following sense.

\begin{defi}[Wedges]\label{def:wedges}
A \emph{wedge} $w$ on a translation surface $X$ is a pair of saddle connections $w=(w_l,w_r)$ such that:
\begin{enumerate}
\item $w_l$ and $w_r$ start from the same conical singularity of $X$;
\item $w_l$ is left-slanted (i.e.\ $\Re(w_l)<0$) and $w_r$ is right-slanted (i.e.\ $\Re(w_r)>0$);
\item $(w_l,w_r)$ consist of two edges of an embedded triangle in $X$.
\end{enumerate}
\end{defi}

A \emph{quadrilateral} $q$ in $X$ is the image of an isometrically embedded quadrilateral in $\CC$ so that the vertices are singularities of $X$, and $q$ contains no other singularities.

\begin{defi}[Admissible quadrangulation]\label{def:quadrangulation}
A quadrilateral $q$ in $X$ is \emph{admissible} if left-slanted and right-slanted saddle connections alternate while we turn around the quadrilateral.

A \emph{quadrangulation} $Q$ of $X$ is a decomposition of $X$ into a union of admissible quadrilaterals.
\end{defi}

Given a quadrilateral $q\in Q$, let us call the saddle connections that start from the same singularity the bottom sides of $q$ and the ones that end on the same singularities the top sides.
We remark that the bottom sides of an admissible quadrilateral, such as one in a quadrangulation of $X$, form a wedge in the sense of the above definition, which we will call the \emph{base wedge} of $q$.

Let $q$ be an admissible quadrilateral and $w=(w_l,w_r)$ its base wedge.
We say that a $q$ is \emph{left-slanted} if its diagonal is left-slanted.
Equivalently, the outgoing vertical separatrix contained in the base wedge of the quadrilateral crosses the top left side.
Similarly, we say that a $q$ is \emph{right-slanted} if its diagonal is right-slanted.

A \emph{diagonal change} in an admissible quadrilateral $q$ consists in replacing the base wedge $w$ with a new one.
More precisely, if $q$ is left-slanted the new base wedge will be $w'=(w_l,w_d)$, where $w_d$ is the diagonal of $q$ itself.
Similarly, if $q$ is right-slanted, the new base wedge will be $w'=(w_d,w_r)$.
Remark that in both cases, thanks to our assumption on the slantedness of $q$ the new base wedge still contains a vertical outgoing separatrix.
To coherently combine diagonal changes in different quadrilaterals, we will need one more geometrical definition.

\begin{defi}[Staircases]\label{def:staircases}
Given a quadrangulation $Q$ of $X$ a \emph{left staircase $S$ for $Q$} (respectively a \emph{right staircase $S$ for $Q$}) is a subset $S\subset X$ which is the union of quadrilaterals $q_1,\dots,q_n$ of $Q$ that are cyclically glued so that the top left (resp.\ top right) side of $q_i$ is identified with the bottom right (resp.\ bottom left) side of $q_{i+1}$ for $1\leq i<n$ and of $q_1$ for $i=n$.

A left (respectively right) staircase $S$ is \emph{well slanted} if all its quadrilaterals are left (resp.\ right) slanted.
\end{defi}

\begin{defi}[Staircase move]\label{def:staircasemove}
Given a quadrangulation $Q$ and a well-slanted left staircase (respectively a well-slanted right staircase) $S$, the \emph{staircase move} in $X$ is the operation which consists in doing simultaneously left (resp.\ right) diagonal changes in all the quadrilaterals of $X$.
\end{defi}

Having given the basic definitions of the diagonal changes algorithm, we now proceed describing the formalism used to encode it.

\begin{defi}[Combinatorial datum]\label{def:combinatorial}
Let $Q$ be a quadrangulation of $k$ quadrilaterals.
Let $q_i$ denote the quadrilateral labeled by $i\in \{ 1,\dots,k\}$.
The \emph{combinatorial datum} $\underline{\pi}=\underline{\pi}_Q$ of the labeled quadrangulation $Q$ is a pair $(\pi_l,\pi_r)$ of permutations of $\{1,\dots,k\}$ such that:
\begin{enumerate}
\item for each $1\leq i \leq k$, the top left side of $q_i$ is glued with the bottom right side of $q_{\pi_l(i)}$;
\item for each $1\leq i \leq k$, the top right side of $q_i$ is glued with the bottom left side of $q_{\pi_r(i)}$;
\end{enumerate}
\end{defi}

We remark that, since $w_{i,l}$ and $w_{\pi_l(i),r}$ are the left sides of the quadrilateral $q_i$ and $w_{i,r}$ and $w_{\pi_r(i),l}$ are its right sides, we have
\begin{equation}\label{eq:traintrackrel}
    w_{i,l}+w_{\pi_l(i),r}=w_{i,r}+w_{\pi_r(i),l},     \qquad \text{ for } 1\leq i \leq k.
\end{equation}
These equations are called \emph{train-track relations}.

Conversely, we can construct a surface with an admissible quadrangulation, starting with a pair of permutations of $k$ elements $\underline{\pi}=(\pi_l,\pi_r)$ and a length datum
\[
    \underline{w}=((w_{1,l},w_{1,r}),\dots,(w_{k,l},w_{k,r}))\in ((\RR_-\times\RR_+)\times(\RR_+\times\RR_+))^k,
\]
where $\RR_-=\set{t\in\RR: t<0}$ and $\RR_+=\set{t\in\RR: t\geq 0}$.
If $\underline{w}$ satisfies the train-track relations~\eqref{eq:traintrackrel} we can build a labeled quadrangulation $Q$ that we denote $(\underline{\pi},\underline{w})$.

We remark that in~\cite{DelecroixUlcigrai:diagonalchanges}, horizontal vectors are not allowed in a quadrangulation, as this would not allow the definition of a \emph{backward} diagonal changes algorithm
However, since we will only use the algorithm forward it will be useful to allow for horizontal saddle connections in the quadrangulation.

\subsection{Moves and matrices}
Let $Q=(\underline{\pi},\underline{w})$ be a labeled quadrangulation.
For each quadrilateral $q_i\in Q$, let $(w_{i,l},w_{i,r})$ be its base wedge and call $w_d$ its diagonal, given by
\[
    w_{i,d} = w_{i,l}+w_{\pi_l(i),r}=w_{i,r}+w_{\pi_r(i),l},
\]
where the equality holds thanks to~\eqref{eq:traintrackrel}.
Given a cycle $c$ of a permutation $\pi_r$ the corresponding staircase $S_c$ formed by the quadrilaterals labeled by the elements of $c$ is well-slanted only if $\Re(w_{i,d})<0$ for all $i\in c$, and similarly if $c$ is a cycle of $\pi_l$.

Starting from a cycle $c$ of $\pi_r$, if its staircase $S_c$ is well-slanted, we can perform a staircase diagonal change as in \cref{def:staircasemove}.
The new length data $\underline{w}'$ is given by
\[
    w_i'=\begin{cases}
            (w_{i,d},w_{i,r}),    &    \text{if $i\in c$;}\\
             w_i,                 &    \text{otherwise}.
         \end{cases}
\]
The new combinatioral datum $\underline{\pi}'=(\pi_l',\pi_r')$ of the new quadrangulation $Q'$ is given by
\begin{equation}\label{eq:piprimeleft}
    \pi_l'(i)=\begin{cases}
                \pi_l\circ\pi_r(i),    &    \text{if $i\in c$;}\\
                \pi_l(i),              &    \text{otherwise.}
              \end{cases}
    \qquad
    \text{and}
    \qquad
    \pi_r'=\pi_r.
\end{equation}

Similarly, if $c$ is a cycle of $\pi_l$ and the corresponding staircase $S_c$ is well-slanted, the new quadrangulation $Q'=(\underline{\pi}',\underline{w}')$ will be given by
\[
    w_i'=\begin{cases}
            (w_{i,l},w_{i,d}),    &    \text{if $i\in c$;}\\
             w_i,                 &    \text{otherwise}.
         \end{cases}
\]
and
\begin{equation}\label{eq:piprimeright}
    \pi_r'(i)=\begin{cases}
                \pi_r\circ\pi_l(i),    &    \text{if $i\in c$;}\\
                \pi_r(i),              &    \text{otherwise.}
              \end{cases}
    \qquad
    \text{and}
    \qquad
    \pi_l'=\pi_l.
\end{equation}

We remark that the operation on the combinatorial datum does not depend on the length datum and that the operation on the wedges $\underline{w}$ is linear.
Hence we can write $\underline{\pi}'=c\cdot\underline{\pi}$, where the action is described above, and we can introduce matrices to encode the action on the length datum.
These matrices will be denoted by $A_{\underline{\pi},c}\in\SL(2k,\ZZ)$.
Let us index the rows and columns of $A_{\underline{\pi},c}$ with $(1,l),(1,r),\dots,(k,l),(k,r)$.
Denote $I_{2k}$ the $2k\times 2k$ identity matrix and for $1\leq i,j\leq k$, and $\epsilon,\nu\in\{l,r\}$ let $E_{(i,\epsilon),(j,\nu)}$ be the $2k\times 2k$ matrix whose entry in row $(1,\epsilon)$ and column $(j,\nu)$ is $1$ and all the other entries are $0$.
We set
\begin{equation}\label{eq:elementarymatrices}
    A_{\underline{\pi},c}=\begin{cases}
                            I_{2k}+\sum_{i\in c} E_{(i,l),(\pi_l(i),r)},
                                & \text{ if $c$ is a cycle of $\pi_r$};\\
                            I_{2k}+\sum_{i\in c} E_{(i,r),(\pi_r(i),l)},
                                & \text{ if $c$ is a cycle of $\pi_l$}.
                           \end{cases}
\end{equation}

Let us summarize the previous discussion.

\begin{lemma}[Staircase move on data]\label{lemma:data}
Given a labeled quadrangulation $Q=(\underline{\pi},\underline{w})$ and a cycle $c$ of $\underline{\pi}$, if the staircase $S_c$ is well slanted, when performing on $Q$ the staircase move in $S_c$ one obtains a new labeled quadrangulation $Q'=(\underline{\pi}',\underline{w}')$ with
\[
    \underline{\pi}' = c \cdot \underline{\pi}, \qquad
    \underline{w}'=A_{\underline{\pi},c}\underline{w},
\]
where $c\cdot\underline{\pi}$ and $A_{\underline{\pi},c}$ are given by \cref{eq:piprimeleft,eq:piprimeright,eq:elementarymatrices}.
\end{lemma}

\subsection{A simpler description of diagonal changes in $\cH(2)$}
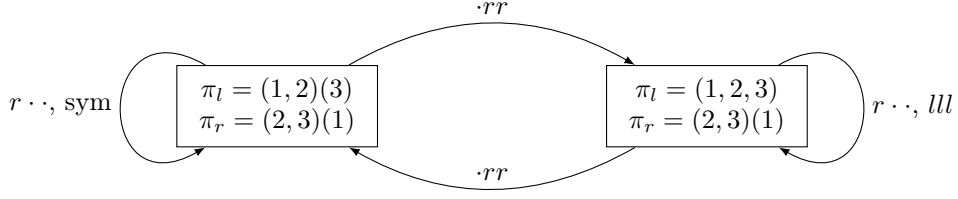
\begin{figure}[bt]
	\begin{tikzpicture}[>=stealth]
     \node (left) [left=1.5cm, shape=rectangle, draw]{
		 \begin{tabular}{c}
			 $\pi_l=(1,2)(3)$ \\
			 $\pi_r=(2,3)(1)$
		 \end{tabular}
     };
     \node (right) [right=1.5cm, shape=rectangle, draw]{
		 \begin{tabular}{c}
			$\pi_l=(1,2,3)$\\
			$\pi_r=(2,3)(1)$
		\end{tabular}
     };
     \draw[->,>=latex] (left) to[bend left]
     node[above] {$\cdot rr$}
     (right);
     \draw[->,>=latex] (right) to[bend left]
     node[above] {$\cdot rr$}
     (left);
     \draw[->,>=latex] (right) to [out=30, in=330, looseness=4]
     node[right]{$r\cdot\cdot$, $lll$}
     (right);
     \draw[->,>=latex] (left) to [out=150, in=210, looseness=4]
     node[left] {$r\cdot\cdot$, sym}
     (left);
   \end{tikzpicture}
	 \caption{The possible moves in $\cH(2)$, up to relabeling and symmetry.}
	 \label{fig:graphH2}
\end{figure}

A more convenient description, for our purposes, of diagonal changes in $\cH(2)$ is given by the following.
Let us introduce a move, called \emph{symmetry}, which exchange the left and right vectors in every quadrilateral.
Moreover, we allow to relabel the wedges.
The graph we obtain is drawn in \cref{fig:graphH2}.

We have introduced these extra moves for the following reasons.
The octagonal continued fraction constructed in~\cite{SmillieUlcigrai:geodesic} and recalled in \cref{sec:octagonFarey} uses also orientation reversing affine diffeomorphisms.
Hence the symmetry is needed in order to represent via diagonal changes that algorithm, precisely for the moves corresponding to even numbered sectors.
Moreover, since the moves of the octagon Farey map act on the unlabeled quadrangulation of the octagon, we need to forego that extra data, that is we have to allow for relabelings.
In fact, the (combinations of) moves of the diagonal changes that correspond to the octagon Farey map, usually begin at one vertex of the graph of possible moves in $\cH(2)$ and end at one which is different from the original one.
The starting vertex and the final one differ precisely by a relabeling.
Relabeling the wedges hence is needed to make sure that the concatenation of diagonal changes agrees with the action of the octagon Farey map; and also allows us to combine the moves from one step to the next.

In the basis given by $\{ E_{(1,l)}, E_{(1,r)},\ldots,E_{(3,r)} \}$, the moves in \cref{fig:graphH2} are given by the following matrices.
\begin{itemize}
	\item $\cdot rr$ from the left node to the right one:
		\[
			\begin{pmatrix}
				1 & 0 & 0 & 0 & 0 & 0 \\
				0 & 1 & 0 & 0 & 0 & 0 \\
				0 & 1 & 1 & 0 & 0 & 0 \\
				0 & 0 & 0 & 1 & 0 & 0 \\
				0 & 0 & 0 & 0 & 1 & 1 \\
				0 & 0 & 0 & 0 & 0 & 1
			\end{pmatrix}
		\]
	\item $\cdot rr$ from the right node to the left one:
		\[
			\begin{pmatrix}
				1 & 0 & 0 & 0 & 0 & 0 \\
				0 & 1 & 0 & 0 & 0 & 0 \\
				0 & 0 & 1 & 0 & 0 & 1 \\
				0 & 0 & 0 & 1 & 0 & 0 \\
				0 & 1 & 0 & 0 & 1 & 0 \\
				0 & 0 & 0 & 0 & 0 & 1
			\end{pmatrix}
		\]
	\item $r\cdot\cdot$ (which is the same matrix in both nodes):
		\[
			\begin{pmatrix}
				1 & 0 & 0 & 1 & 0 & 0 \\
				0 & 1 & 0 & 0 & 0 & 0 \\
				0 & 0 & 1 & 0 & 0 & 0 \\
				0 & 0 & 0 & 1 & 0 & 0 \\
				0 & 0 & 0 & 0 & 1 & 0 \\
				0 & 0 & 0 & 0 & 0 & 1
			\end{pmatrix}
		\]
	\item $lll$ plus relabeling:
		\[
			\begin{pmatrix}
				0 & 0 & 1 & 0 & 0 & 0 \\
				0 & 0 & 0 & 1 & 1 & 0 \\
				0 & 0 & 0 & 0 & 1 & 0 \\
				0 & 0 & 1 & 0 & 0 & 1 \\
				1 & 0 & 0 & 0 & 0 & 0 \\
				1 & 1 & 0 & 0 & 0 & 0
			\end{pmatrix}
		\]
	\item the left/right symmetry plus relabeling:
		\[
			\begin{pmatrix}
			0 & 0 & 0 & 0 & 0 & 1 \\
		 	0 & 0 & 0 & 0 & 1 & 0 \\
			0 & 0 & 0 & 1 & 0 & 0 \\
			0 & 0 & 1 & 0 & 0 & 0 \\
			0 & 1 & 0 & 0 & 0 & 0 \\
			1 & 0 & 0 & 0 & 0 & 0
			\end{pmatrix}
		\]
\end{itemize}

\section{The octagon Farey map in terms of diagonal changes}\label{sec:thm}
In this section, we show that the octagon Farey map $F$ is an acceleration of diagonal changes moves.
Given a direction $\theta=[s_0;s_1,s_2,\ldots]_\cO$, we have a well-defined sequence of maps $(F_{s_i})_{s_i\in\NN}$.
These maps act affinely on the surface $X_\cO$.
As we said above, with a slight abuse of notation, we will refer to this sequence of maps also as the octagon Farey map.

\begin{figure}[t]
\centering
\def\svgwidth{0.9\textwidth}
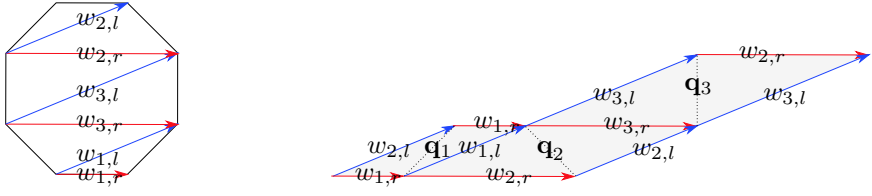
\caption{The quadrangulation $Q_0$ of the regular octagon.}
\label{fig:Q0}
\end{figure}

By definition of the octagon Farey map the first entry $s_0$ plays a special role.
Since $s_0$ determines in which of the eight sectors $\overline{\Sigma}_j$ lies the direction $\theta$, this determines the starting quadrangulation of the surface $X_\cO$.
More precisely, let $Q_0$ be the quadrangulation in \cref{fig:Q0}.
Then the beginning quadrangulation of $X_\cO$ is $Q=\nu_{s_0}^{-1} Q_0$.

\begin{figure}[b]
\centering
\def\svgwidth{0.8\textwidth}
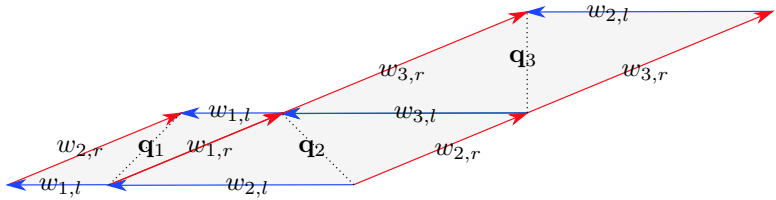
\caption{The beginning quadrangulation $Q'$ of the regular octagon.}
\label{fig:openedwedges}
\end{figure}

We now describe how to translate the induced action of the octagon Farey map in terms of diagonal changes.
We remark that, as $s_0$ dictates the starting quadrangulation, and the other $s_i$ only take values from $1$ to $7$, we only have to translate these seven cases.
In order to exploit the symmetry among the seven sectors $\overline{\Sigma}_j$, $j=1,\ldots, 7$, and in order to make clearer pictures, we apply the map $F_{s_0}$ to $X_\cO$, thus opening up the sector $\overline{\Sigma}_{s_0}$ onto $\overline{\Sigma}$.
The quadrangulation $Q'=\gamma Q_0$ we obtain is as in \cref{fig:openedwedges}.
We label quadrilaterals in $Q'$ so that the combinatorial datum $\underline{\pi}$ is given by
\[
    \pi_l=(1,2)(3), \qquad \text{and} \qquad \pi_r=(1)(2,3).
\]
Let us remark that we will use diagonal changes to approximate the direction $\theta=[s_1; s_2, \ldots]_\cO$ and \emph{not} the vertical one.
Figures that represent the movements can be found at the end of the document, see \cref{sec:figures} for some comments about them.

\begin{enumerate}
		\item First sector  $\bigl(\frac{\pi}{8}\leq\theta\leq\frac{2\pi}{8}\bigr)$: $\cdot rr$, $r\cdot\cdot$, $\cdot rr$, see \cref{fig:firstsector}.
		Hence:
		\[
	    A_1  =   \begin{pmatrix}
	                1 & 0 & 0 & 1 & 0 & 0\\
	                0 & 1 & 0 & 0 & 0 & 0\\
	                0 & 1 & 1 & 0 & 0 & 1\\
	                0 & 0 & 0 & 1 & 0 & 0\\
	                0 & 1 & 0 & 0 & 1 & 1\\
	                0 & 0 & 0 & 0 & 0 & 1
	             \end{pmatrix}.
		\]

		\item Second sector  $\bigl(\frac{2\pi}{8}\leq\theta\leq\frac{3\pi}{8}\bigr)$: $\cdot rr$, $lll$, $r\cdot r$, $\cdot r \cdot$, symmetry, see \cref{fig:secondsector}.
		Hence:
		\[
    	A_2 =   \begin{pmatrix}
                1 & 1 & 0 & 0 & 0 & 0 \\
								1 & 0 & 0 & 1 & 1 & 1 \\
								0 & 1 & 1 & 0 & 0 & 1 \\
								1 & 1 & 0 & 0 & 1 & 1 \\
								0 & 0 & 0 & 1 & 1 & 1 \\
								0 & 2 & 2 & 0 & 0 & 1
             	\end{pmatrix}.
		\]

		\item Third sector  $\bigl(\frac{3\pi}{8}\leq\theta\leq\frac{4\pi}{8}\bigr)$: $\cdot rr$, $lll$, $lll$, $\cdot rr$, see \cref{fig:thirdsector}.
		Hence:
		\[
    	A_3 =   \begin{pmatrix}
                0 & 0 & 0 & 0 & 1 & 1 \\
								1 & 1 & 1 & 0 & 0 & 1 \\
								1 & 1 & 1 & 1 & 1 & 1 \\
                1 & 1 & 0 & 0 & 1 & 1 \\
                1 & 2 & 2 & 0 & 0 & 1 \\
                0 & 1 & 1 & 1 & 1 & 1
             	\end{pmatrix}.
		\]

		\item Fourth sector  $\bigl(\frac{4\pi}{8}\leq\theta\leq\frac{5\pi}{8}\bigr)$: $\cdot\cdot l$, $\cdot rr$, $\cdot rr$, $ll\cdot$, $ll\cdot$, $r\cdot\cdot$, symmetry, see \cref{fig:fourthsector}.

		These moves correspond in the reduced graph to:
		symmetry, $r\cdot\cdot$, symmetry, $\cdot rr$, $\cdot rr$, symmetry $\cdot rr$, $\cdot rr$, symmetry, $r\cdot\cdot$, symmetry.
		Hence:
		\[
    	A_4 =   \begin{pmatrix}
								0 & 0 & 1 & 0 & 0 & 1 \\
								0 & 1 & 1 & 0 & 1 & 1 \\
								1 & 1 & 1 & 1 & 1 & 1 \\
								0 & 1 & 2 & 0 & 0 & 1 \\
								1 & 2 & 1 & 0 & 1 & 1 \\
                2 & 1 & 1 & 1 & 1 & 1
             \end{pmatrix}.
		\]

		\item Fifth sector  $\bigl(\frac{5\pi}{8}\leq\theta\leq\frac{6\pi}{8}\bigr)$: $\cdot\cdot l$, $\cdot rr$, $lll$, $r\cdot r$, $l\cdot \cdot$, see \cref{fig:fifthsector}.

		These moves correspond in the reduced graph to:
		symmetry, $r\cdot\cdot$, symmetry, $\cdot rr$, $lll$, $\cdot rr$, symmetry, $r\cdot\cdot$, symmetry.
		Hence:
		\[
    	A_5 =   \begin{pmatrix}
								0 & 1 & 1 & 0 & 0 & 0 \\
								0 & 0 & 1 & 1 & 1 & 1 \\
								1 & 1 & 1 & 0 & 1 & 1 \\
								0 & 1 & 2 & 0 & 0 & 1 \\
								1 & 0 & 1 & 1 & 1 & 1 \\
								2 & 2 & 1 & 0 & 1 & 1 \\
             \end{pmatrix}.
		\]

		\item Sixth sector  $\bigl(\frac{6\pi}{8}\leq\theta\leq\frac{7\pi}{8}\bigr)$: $ll\cdot$, $\cdot\cdot l$, $rrr$, $l\cdot l$, see \cref{fig:sixthsector}.

		These moves correspond in the reduced graph to: symmetry, $\cdot rr$, $r\cdot\cdot$, $lll$, $\cdot rr$.
		Hence:
		\[
    	A_6 =   \begin{pmatrix}
								0 & 0 & 0 & 1 & 1 & 0 \\
								1 & 1 & 1 & 0 & 0 & 0 \\
								1 & 1 & 1 & 0 & 1 & 1 \\
								1 & 0 & 0 & 1 & 1 & 0 \\
                1 & 1 & 2 & 0 & 0 & 1 \\
                0 & 0 & 1 & 0 & 1 & 1
             \end{pmatrix}.
		\]

		\item Seventh sector $\bigl(\frac{7\pi}{8}\leq\theta\leq\pi\bigr)$: $ll\cdot$, $\cdot\cdot l$, $ll\cdot$, $\cdot\cdot l$, see \cref{fig:seventhsector}.

		These moves correspond in the reduced graph to: symmetry, $\cdot rr$, $r\cdot\cdot$, $\cdot rr$, $r\cdot\cdot$, symmetry.
		Hence:
		\[
    	A_7 =   \begin{pmatrix}
                1 & 0 & 0 & 0 & 0 & 0\\
                1 & 1 & 0 & 0 & 1 & 0\\
                0 & 0 & 1 & 0 & 0 & 0\\
                1 & 0 & 0 & 1 & 1 & 0\\
                0 & 0 & 0 & 0 & 1 & 0\\
                0 & 0 & 2 & 0 & 0 & 1
             \end{pmatrix}.
		\]
\end{enumerate}

\section*{Acknowledgments}
The author thanks Corinna Ulcigrai and Vincent Delecroix for useful conversation.
I would also like to thank Davide Ravotti for comments on an earlier draft.

\printbibliography

\appendix
\section{Drawings}\label{sec:figures}
In the last few pages of this document we present the drawings that describe the concatenation of diagonal changes moves needed to recover the octagon Farey map.

Let us comment on the pictures that follows.
In every picture we represent at the top the quadrangulation $Q'$ together with a line in a generic direction $\theta$ inside the appropriate sector.
Then we represent the diagonal changes in left to right, top to bottom order.
In order to keep the pictures as clear as possible, labels are kept to a minimum and we do not represent the direction $\theta$ in the drawings of the staircases moves.
The reader can check that all the moves are admisibles, that is, the staircases are slanted in the appropriate direction.
Moreover, in order to save space, we do not represent the final symmetry move in the even numbered sectors.

\begin{figure}[bt]
	\centering
	\def\svgwidth{0.8\textwidth}
	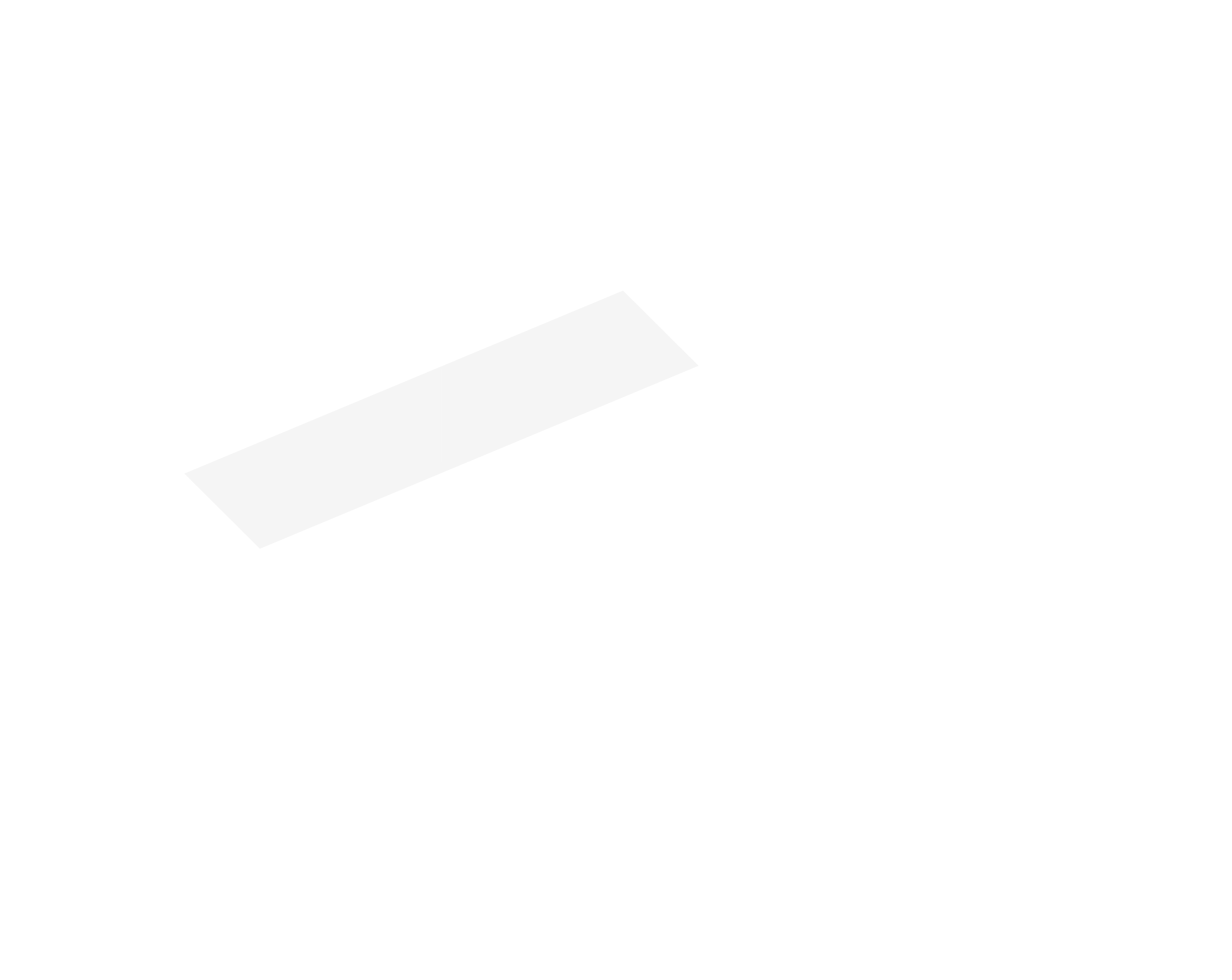
	\caption{The moves of the diagonal changes algorithm for the first sector.}
	\label{fig:firstsector}
	\def\svgwidth{0.8\textwidth}
	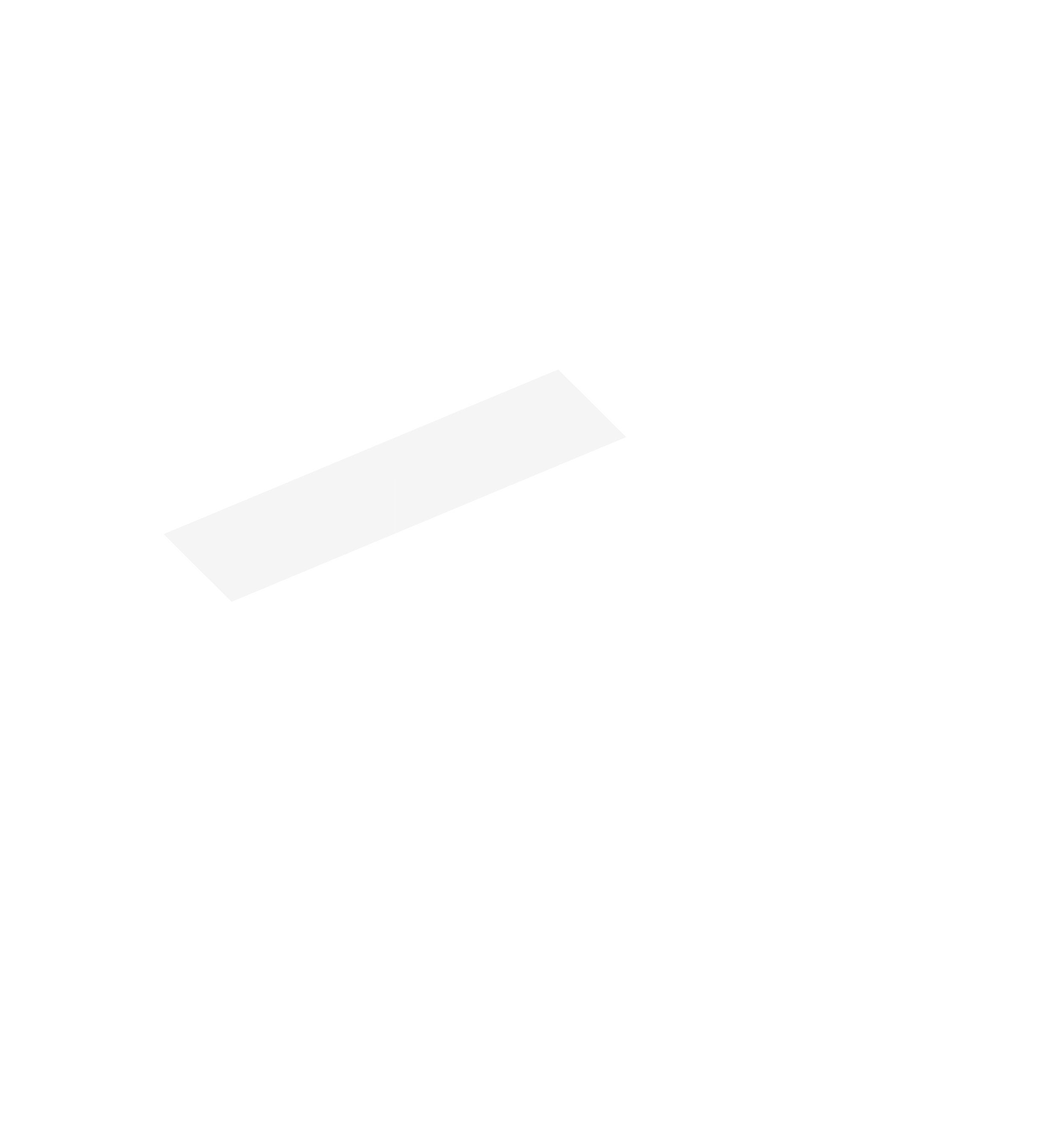
	\caption{The moves of the diagonal changes algorithm for the second sector.}
	\label{fig:secondsector}
\end{figure}

\begin{figure}[bt]
	\centering
	\def\svgwidth{0.8\textwidth}
	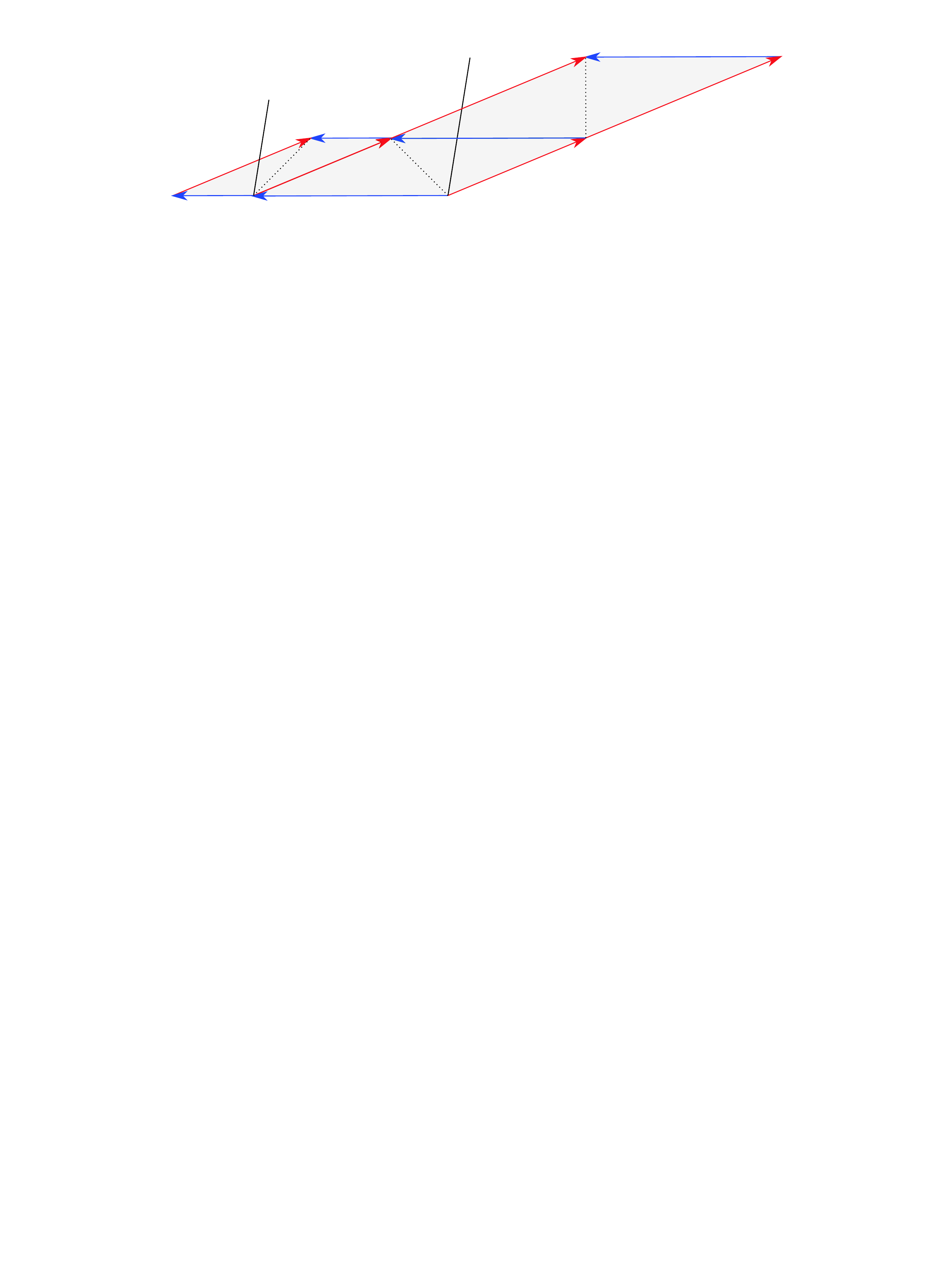
	\caption{The moves of the diagonal changes algorithm for the third sector.}
	\label{fig:thirdsector}
\end{figure}

\begin{figure}[bt]
	\centering
	\def\svgwidth{0.8\textwidth}
	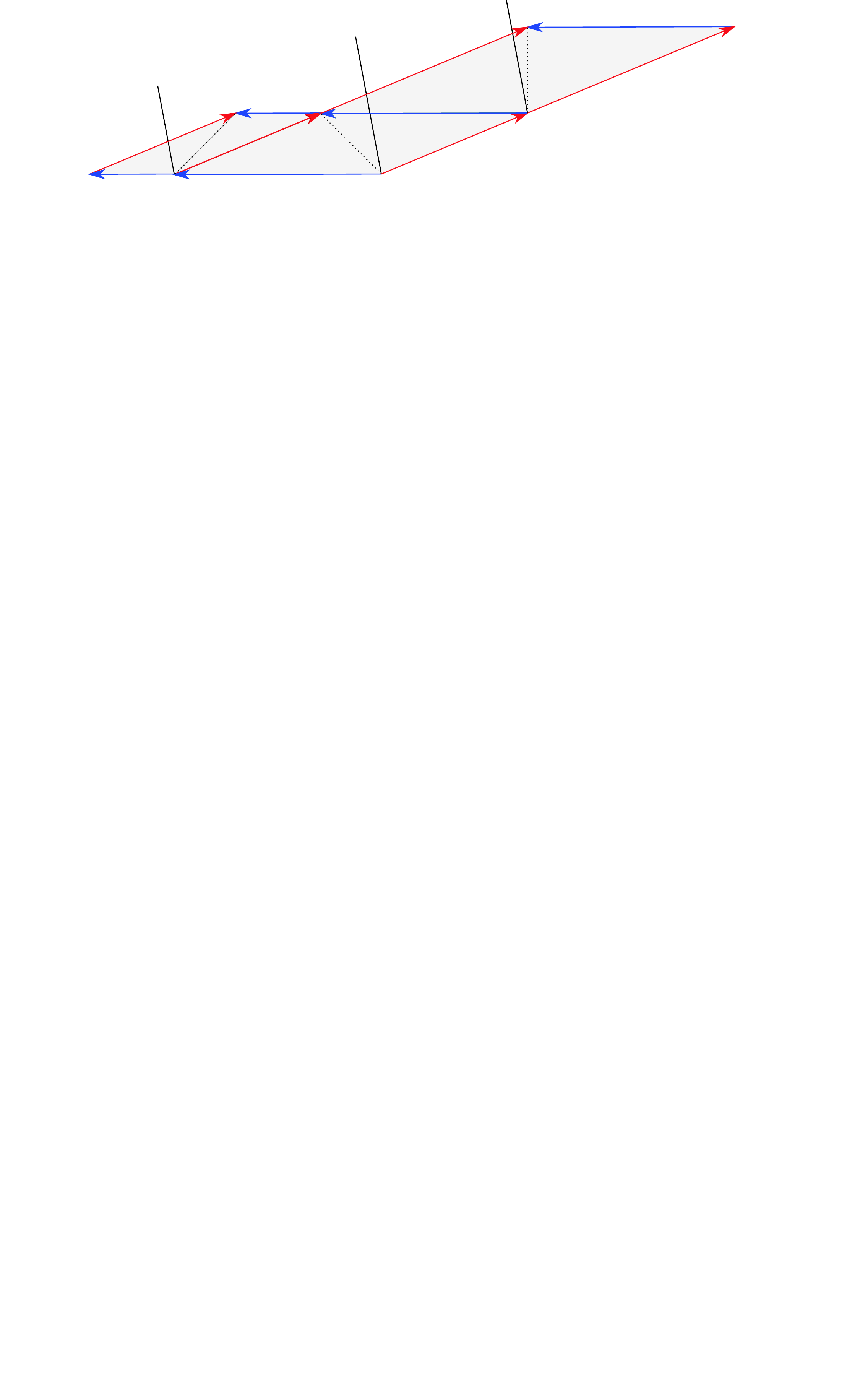
	\caption{The moves of the diagonal changes algorithm for the fourth sector.}
	\label{fig:fourthsector}
\end{figure}

\begin{figure}[bt]
	\centering
	\def\svgwidth{0.8\textwidth}
	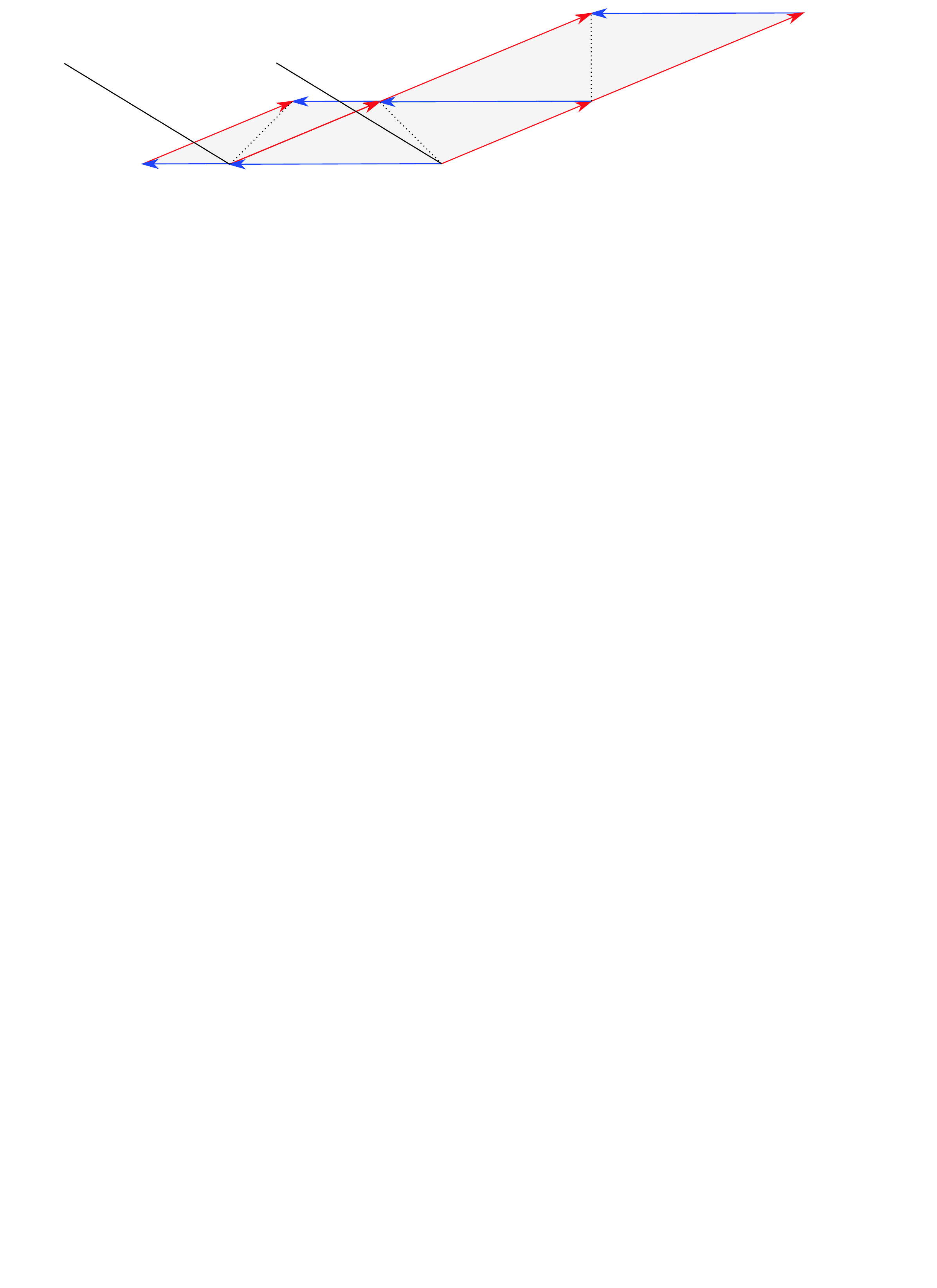
	\caption{The moves of the diagonal changes algorithm for the fifth sector.}
	\label{fig:fifthsector}
\end{figure}

\begin{figure}[bt]
	\centering
	\def\svgwidth{0.8\textwidth}
	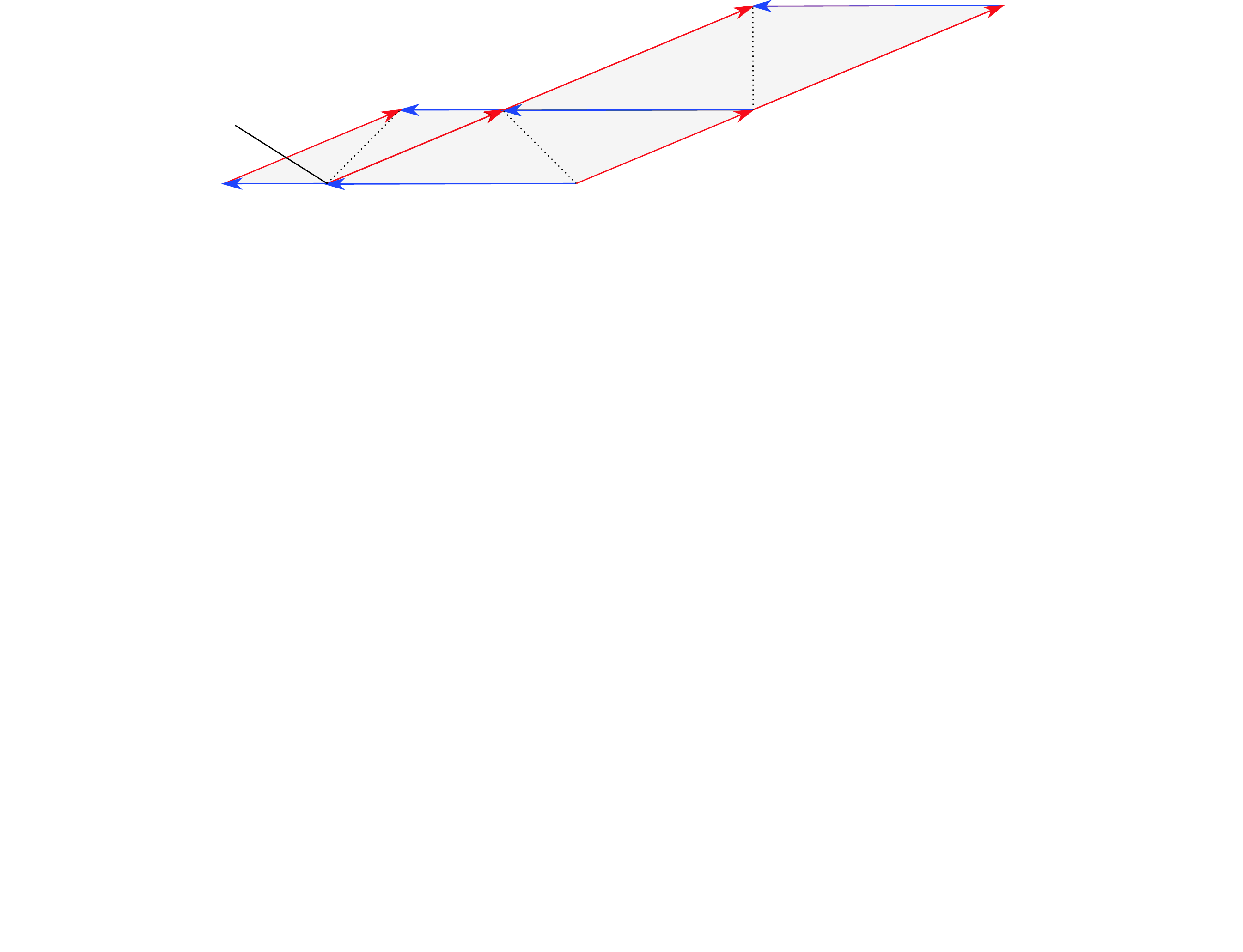
	\caption{The moves of the diagonal changes algorithm for the sixth sector.}
	\label{fig:sixthsector}
	\def\svgwidth{0.8\textwidth}
	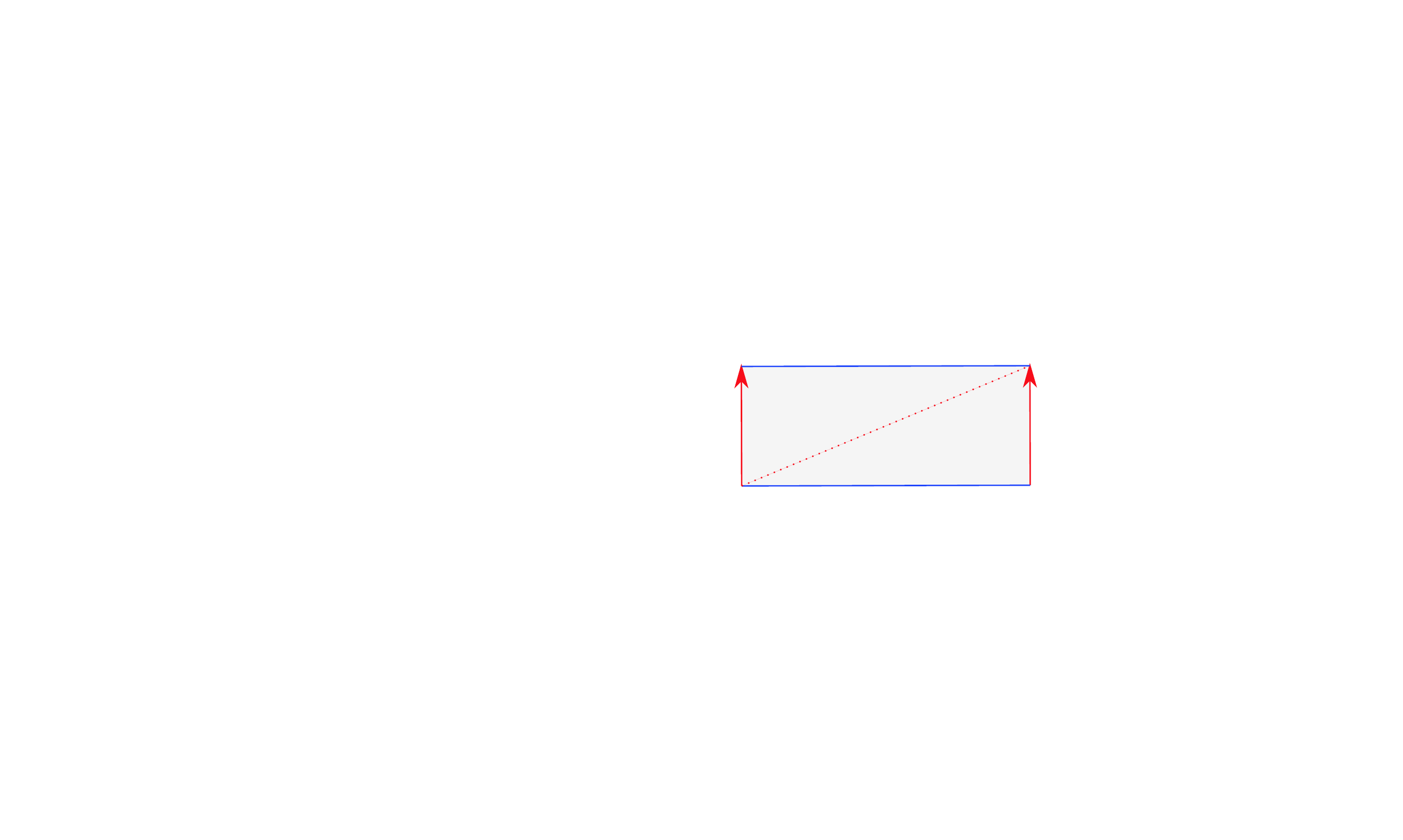
	\caption{The moves of the diagonal changes algorithm for the seventh sector.}
	\label{fig:seventhsector}
\end{figure}

\end{document}